\documentclass[preprint]{elsarticle}

\usepackage{amssymb}
\usepackage{geometry}
\usepackage{lmodern}
\usepackage{amsmath,amssymb}

\usepackage[T1]{fontenc}
\usepackage[utf8]{inputenc}
\usepackage{appendix}
\usepackage{paralist}
\usepackage{longtable}
\usepackage{graphicx}
\usepackage{amsmath}
\usepackage{fullwidth}
\usepackage{xcolor}

\usepackage{amssymb}
\usepackage{xcolor}
\usepackage{upgreek}
\usepackage{marvosym}
\usepackage{amsthm}
\usepackage{caption}
\usepackage{epsfig,epstopdf}
\usepackage{subcaption}
\usepackage{amstext}
\usepackage{array}
\usepackage{euler}
\usepackage{times}
\usepackage{lipsum}
\usepackage{xspace}

\numberwithin{equation}{section}
\usepackage{algorithm}% http://ctan.org/pkg/algorithms
\usepackage{algorithmic}% http://ctan.org/pkg/algorithms

\numberwithin{algorithm}{section}
\theoremstyle{plain}% Theorem-like structures provided by amsthm.sty
\newtheorem{theorem}{Theorem}[section]
\newtheorem{lemma}[theorem]{Lemma}
\newtheorem{corollary}[theorem]{Corollary}
\newtheorem{proposition}[theorem]{Proposition}

\theoremstyle{definition}
\newtheorem{definition}[theorem]{Definition}

\usepackage{amsmath, amssymb} % Make sure this is in your preamble
\newcommand{\R}{\mathbb{R}}   % Shortcut for the real numbers
\newcommand{\calH}{\mathcal{H}} % Shortcut for calligraphic H

\theoremstyle{remark}
\newtheorem{remark}{Remark}

\usepackage{color}
\journal{}

\begin{document}

\begin{frontmatter}

%% Title, authors and addresses

%% use the tnoteref command within \title for footnotes;
%% use the tnotetext command for theassociated footnote;
%% use the fnref command within \author or \address for footnotes;
%% use the fntext command for theassociated footnote;
%% use the corref command within \author for corresponding author footnotes;
%% use the cortext command for theassociated footnote;
%% use the ead command for the email address,
%% and the form \ead[url] for the home page:
%% \title{Title\tnoteref{label1}}
%% \tnotetext[label1]{}
%% \author{Name\corref{cor1}\fnref{label2}}
%% \ead{email address}
%% \ead[url]{home page}
%% \fntext[label2]{}
%% \cortext[cor1]{}
%% \affiliation{organization={},
%%             addressline={},
%%             city={},
%%             postcode={},
%%             state={},
%%             country={}}
%% \fntext[label3]{}

\title{A Decomposition Method for Finite-Time Stabilization of Bilinear  Systems with Applications to Parabolic and Hyperbolic Equations}

\author[1]{Kamal FENZA}
\author[2]{Moussa LABBADI}
\author[1]{Mohamed OUZAHRA}

\address[1]{USMBA - Ecole Normale Supérieure (ENS), Fez, 30010, Morocco}
\address[2]{Aix-Marseille University, LIS, UMR CNRS 7020, Marseille, 13013, France}
%\author[a]{ Kamal FENZA}
%\author[b]{Moussa LABBADI}
%\author[a]{Mohamed OUZAHRA}
%
%
%
%
%\affiliation[a]{organization={USMBA - Ecole Normale Supérieure (ENS)}, 
%                city={Fez},
%                postcode={30010},
%                country={Morocco}}
%
%\affiliation[b]{organization={Aix-Marseille University, LIS, UMR CNRS 7020}, 
%                city={Marseille},
%                postcode={13013},
%                country={France}}
                
\begin{abstract}
In this work, we address the problem of finite-time stabilization for a class of  bilinear system. We propose a decomposition-based approach in which the nominal system is split into two subsystems, one of which is inherently finite-time stable without control. This allows the stabilization analysis to focus solely on the remaining subsystem. To ensure the well-posedness of the closed-loop system, we establish sufficient conditions on the system and control operators. The stabilization results are then derived using a suitable Lyapunov function and an observation condition. The effectiveness of the proposed approach is demonstrated through examples involving both  parabolic and hyperbolic infinite-dimensional systems.
\end{abstract}

% %%Graphical abstract
% \begin{graphicalabstract}
% %\includegraphics{grabs}
% \end{graphicalabstract}

% %%Research highlights
% \begin{highlights}

% \item Well posedness is ensured using semilinear Cauchy problem theory.
% \item Finite time stability is estabished using Lyapunov function.
% \item Applications to both linear and nonlinear systems are treated.
% \end{highlights}

\begin{keyword}
%% keywords here, in the form: keyword \sep keyword
Bilinear systems \sep  Linear systems \sep  Finite-time stabilization  \sep Nonlinear control \sep Space decomposition.
%% PACS codes here, in the form: \PACS code \sep code

%% MSC codes here, in the form: \MSC code \sep code
%% or \MSC[2008] code \sep code (2000 is the default)

\end{keyword}

\end{frontmatter}

%% \linenumbers

%% main text
\section{Introduction}
In this paper, we consider the following  bilinear  system:
\begin{eqnarray}\label{systSvvb}
	\begin{cases}
		\frac{d }{d t} y(t)=A y(t)+ u(t)By(t) , \quad    \hspace*{0.1cm} t>0,\\
		y(0)=y_{0},    \\
	\end{cases}
\end{eqnarray}
where  $y(.)\in \calH$ is the  state  in Hilbert space   and the control input  $u \in \R$ is  a scalar function. The    unbounded linear operator  $A:D(A)\rightarrow \mathcal{H}$ is an infinitesimal generator  of  a  strongly continuous semigroup $S(t)$  on $\mathcal{H}$ and $B:\mathcal{H}\rightarrow \mathcal{H}$ is a linear bounded  operator. 

Bilinear  system~\eqref{systSvvb} can be modeled several real applications in economics, biology and physics (see \cite{APP1, APP2}). For these reasons, this class of systems has received lots of consideration the last thirty years, where good progress is obtained in the asymptotic stabilization  area by various feedback laws. In this paper,   we'll focus on the problem of finite-time stabilization of (\ref{systSvvb}), that consists of the existence of a feedback control and a finite (settling) time $T^{*}$  such that all trajectories of the  (\ref{systSvvb}) converge to zero within this time and remain at zero afterwards. Finite-time stability seems more appropriate for applications in control theory, robotics and stabilization of underwater vehicles, satellites and auto-mobiles, which are subject to severe time constraints. Actually, in the finite-dimensional case, this problem was treated for the  continuous autonomous systems by \cite{Am1}. In addition, they proved the first converse  theorem for the stability in finite time under  the continuity of  the settling time function. In \cite{7}, it was proved that  finite-time stability of  nonlinear systems can be obtained under the homogeneity condition with a negative degree combined with the asymptotic stability. In addition to this, the authors in \cite{14}  studied the finite-time stabilization of   continuous  non autonomous systems  and differential inclusions.  Furthermore,  this problem  has been investigated for the linear systems by \cite{6666666} using  a time-varying feedback.  The problem of  finite-time partial stability  has been investigated in \cite{25}  with a  continuous or discontinuous state feedback laws. In the context of infinite-dimensional systems, the finite-time stability  of  many well-known partial differential equations  is still the topic of intensive research. For example, Polyakov et al (see e.g., \cite{10,15,17})  have extended the notion of homogeneity to infinite-dimensional linear systems with  a positive definite control operator, were,  the authors developed a nonlinear stabilizing feedback and a Lyapunov function for homogeneous dilation to obtain a finite-time stability result, which is then applied to the wave and the heat equations.

  For the case of bilinear reaction-diffusion equation with the control operator $B=I$ (i.e. $I$ the identity of $\mathcal{H}$) was studied in  \cite{2333}  by a direct way in a settling time which depends on the initial state. In \cite{6jamazi, sob}, the finite-time stabilization of abstract bilinear control systems with bounded control operators is investigated. These approaches rely on strong assumptions regarding the system operators. For instance, in \cite{6jamazi}, the operator $A$ is assumed to be homogeneous of a given degree $\beta$ while in \cite{sob}, the control operator $B$ is required to be coercive. As a result, the applicability of these results is restricted to a narrow class of systems. Furthermore, \cite{ha} explored finite-time stabilization
  for unbounded bilinear systems, relying on Lyapunov-based methods and assuming  the  control operators is coercive  on the state space $\mathcal{H}$. Recently, using a state-space decomposition by $\ker(B)$ and its orthogonal, \cite{Amaliki} studied the  finite-time stability of parabolic systems with a  unbounded control operator that is coercive only on $\ker(B)$. Additionally, this decomposition method may not be applicable in many situations, as the invariance of  the space $\ker(B)$ under the semigroup $S(t)$ is not easy to guarantee.  
  
Based on the above discussion, existing results on finite-time stability are mainly applicable to distributed  bilinear systems where the control operator is  coercive. In this paper, we focus on a class of bilinear  parabolic and hyperbolic  systems in infinite-dimensional spaces, which can be  observable or not, and   where the control operator $B$ is  not  coercive in any subspace of $\mathcal{H}$, thereby expanding the understanding of finite-time stability to include a broader range of scenarios. Our new approach   consists of decomposing the dynamical system into two subsystems; one of them is finite time stable (without control), so that the study can concentrate on the other. Moreover, the derived results are applied to ensure the finite time stability of a class of linear systems governed by a positive control operator  but not necessarily definite. We rely  on the construction of the Lyapunov function candidate $V$  and  then give   sufficient conditions on the system's parameters to ensure  the following differential inequality:  
 $$
\dot{V}+c V^\theta \leq 0,  \hspace*{0.1cm} \theta \in(0,1), \hspace*{0.1cm}   c>0.
$$
The plan of the paper is as follows: the Section 2 introduces a review on maximal monotone operator and  the well-posedness of some evolution equations.  In Section 3, we start by decomposing the nominal system into two subsystems, the first of which is uncontrolled, while for the second one,  the problem of global finite-time stabilization is achieved  under nonlinear state feedback control. In Section 4 we show that the results obtained can  be applied to linear systems. Section 5  provides   applications to heat, wave, coupled transport-heat and beams equations. The conclusion is given in Section 6.

%%%%%%%%%%%%%%%%%%%%%%%%%%%%%%%%%%%%%%%%%%%%%%%%%%%%%%%%%%%%%%%%%%%%%%%%%%%%%%%%%%%%%%%%%%%%%%%%%%%%%%%%

\section{ Review on maximal monotone operator and  evolution equations}
This section is devoted to some essential properties of maximal monotone operators and  their  relation with  existence and uniqueness theory for  evolution equations.\\
Throughout the paper, $\mathcal{H}$ is a real Hilbert space, whose norm and scalar product are denoted respectively by $ \left\langle . , .\right\rangle $ and $ \lVert .\rVert$. The norm of  linear bounded  operators on $\mathcal{H}$ will also be denoted by $ \lVert .\rVert$ for convenience and if there is no ambiguity.\\ 
 \begin{definition} ([\cite{sobb}, pp. 20-22]):  The operator $\mathcal{A}: D(\mathcal{A}) \subset \mathcal{H} \rightarrow \mathcal{H}$ is said to be monotone if  $\langle \mathcal{A}(x)-\mathcal{A}(y), x-y\rangle \geq 0$, for all $x, y \in D(\mathcal{A})$ (i.e. $-\mathcal{A}$ is dissipative). If in addition,   the graph of $\mathcal{A} $  is not properly contained in the graph of any other monotone operator in $\mathcal{H}$, we say that  $\mathcal{A} $ is maximal monotone.
\end{definition} 
 \begin{theorem} ([\cite{sobb}, p. 34]):\label{theork1}
Let $\mathcal{A}$ be a maximal monotone operator on $\mathcal{H}$ and $\mathcal{B}$ a Lipschitzian monotone operator from $\mathcal{H}$ into $\mathcal{H}$. Then the operator $\mathcal{A}+\mathcal{B}$ is maximal monotone.
\end{theorem} 
\begin{theorem}  ([\cite{sobb}, p. 36]): \label{theork2}
Let $\mathcal{A}$ and $\mathcal{B}$ be two maximal monotone operators on $\mathcal{H}$. If $D(\mathcal{A}) \cap \operatorname{Int}(D(\mathcal{B})) \neq \emptyset$, where $\operatorname{int}(D(\mathcal{B}))$ denotes the interior of $D(\mathcal{B})$, then $\mathcal{A}+\mathcal{B}$ is a maximal monotone operator.
\end{theorem} 
Now, we introduce the notions of weak and strong  solutions of the following evolution equation:
\begin{eqnarray}\label{equation}
	\begin{cases}
		\frac{d}{d t}  x(t)=\mathcal{A} x(t)+f(t) , \quad    \hspace*{0.1cm} t>0,\\
		x(0)=x_{0}\in \mathcal{H},  \\
	\end{cases}
\end{eqnarray}
where $\mathcal{A}: D(\mathcal{A}) \subset \mathcal{H} \rightarrow \mathcal{H}$ is a  (possibly nonlinear) operator and $f$ is a given function.
\\
\begin{definition} ([\cite{sobb}, p. 64]):
 Let $f \in L^1(0, T ; \mathcal{H})$ with $T>0$.  The function $x \in C([0, T], \mathcal{H})$ 
 is called a strong  solution of  (\ref{equation}) if  it satisfies \\
$\bullet$ $x$ is absolutely continuous on any compact subset from $(0, T)$,\\
$\bullet$ $x(t) \in D(\mathcal{A})$ and $\frac{d x}{d t}(t)+\mathcal{A} x(t)= f(t)$, a.e. on $(0, T)$.
\end{definition} 
\begin{definition} ([\cite{sobb}, p. 64]):
A function $x\in C([0, T], \mathcal{H})$ is a weak solution for  (\ref{equation}) if there exist sequences $(f_{n})\in L^1(0, T ; \mathcal{H})$ such that $f_{n}\rightarrow f$ in $ L^1(0, T ; \mathcal{H})$ 
and $(x_{n})\in C([0, T], \mathcal{H})$ such that  $x_{n}$ is a strong
solution  of $\frac{d }{d t} x_{n}(t)=\mathcal{A} x_{n}(t)+f_{n}(t) $ and $x_{n}\rightarrow x$ uniformly on $[0, T]$.
\end{definition} 
The next result establishes the existence and uniqueness of the weak solution of system (\ref{equation}).
\begin{theorem}  ([\cite{sobb}, p.  65]):\label{4}
Let  $\mathcal{A}: D(\mathcal{A}) \subset \mathcal{H} \rightarrow \mathcal{H}$ be a  maximal monotone operator and $f\in L^1(0, T ; \mathcal{H})$. Then, for all $x_{0}\in \overline{D(\mathcal{A})}$  the  system  (\ref{equation})  admits  a unique weak solution $x\in C([0,+\infty), \mathcal{H})$.
\end{theorem} 
The following lemma gives a comparison principle.
\begin{lemma} ([\cite{k}, p. 102])\label{4nb}:
Let $J\subset\mathbb{R}$ and $f:\mathbb{R}\times J\rightarrow \mathbb{R}$ is continuous in $t\geq 0$  and locally Lipschitz  in $J$. For every $t_0\geq 0$, we consider the following  scalar differential equation:
  \begin{equation}\label{2b}
  \frac{d }{d t} x(t)= f(t,x(t)),  \hspace*{0.5cm} x(t_0)=x_0.
  \end{equation}
   Let $\left[t_0, T\right)$  ($t_0< T\leq \infty $) be the maximal interval of existence of the solution $x(t)$ of  (\ref{2b}) and suppose that $x(t) \in J$ for all $t \in\left[t_0, T\right)$. Let    $v(t)$ be a continuous function whose upper right-hand derivative $D^{+} v(t)$ satisfies the differential inequality
$$
D^{+} v(t) \leq f(t, v(t)), \quad v\left(t_0\right) \leq x_0,
$$
with $v(t) \in J$ for all $t \in\left[t_0, T\right)$. Then, $$v(t) \leq x(t),\hspace*{0.1cm} \forall t \in\left[t_0, T\right).$$
\end{lemma}
%%%%%%%%%%%%%%%%%%%%%%%%%%%%%%%%%%%%%%%%%%%%%%%%%%%%%%%%%%%%%%%%%%%%%%%%%%%%%%%%%%%%%%%%%%%%%%%%%%%%%%%%%%%%%%%%%%
\section{ Finite-time stability of bilinear systems}
\subsection{Problem statement}
We consider the following system:
\begin{eqnarray}\label{systS}
	\begin{cases}
		\frac{d }{d t} y(t)=A y(t)+ u(t)By(t) , \quad   \hspace*{0.1cm} t>0,\\
		y(0)=y_{0},   \\
	\end{cases}
\end{eqnarray}
where\\
$\bullet$ $A:D(A)\rightarrow \mathcal{H}$ is the infinitesimal generator of a linear   $C_{0}$-semigroup $S(t)$ on $\mathcal{H}$.\\
$\bullet$ $B:\mathcal{H}\rightarrow \mathcal{H}$ is a bounded linear self-adjoint operator.\\

We consider the system~\eqref{systS}, which the well-posedness of the system at hand will be studied in closed-loop, so before turning our attention to the existence and uniqueness of the solution, we will discuss first the choice of a feedback law.
To have an idea about the feedback law and the Lyapunov function candidate, we formally differentiate the time rate of change of the “energy”
\[
\frac{1}{2} \frac{d}{dt} \|y(t)\|^2 = \langle Ay(t), y(t) \rangle +  u(t) \langle By(t), y(t) \rangle,
\]
Thus, for instance when \( A \) is dissipative,  the dissipativity of the system can be ensured by selecting the control as $ u(t) = - \langle By(t), y(t) \rangle^{\beta}  $ for a suitable  real exponent \( \beta \). Based on related studies and specific examples found  
in the  literature (see \cite{2333} and \cite{17}), a natural and effective  choice  \( \beta \) is \( \beta = -\mu \) with \( 0 < \mu < \frac{1}{2} \). This   choice yields  the following dissipation inequality:
$$
\frac{d}{dt} \|y(t)\|^2 \leq -2 \langle By(t), y(t) \rangle^{1 - \mu}.
$$

Accordingly, two Lyapunov candidate functions arise:
$$
V(y) = \|y\|^2 \quad \text{and} \quad V(y) = \langle By, y \rangle.
$$

The first candidate can be considered in the case of coercive control operator \( B \), this case has been considered in
\cite{sob} and \cite{6jamazi}. Otherwise we consider the second Lyapunov function. Accordingly, the second choice appears more suitable, as it can still be applied when $B$ is  coercive. Taking into account that, in the latter, the	positive function $V$ is not definite, one can not immediately
conclude about the finite-time convergence of the state, other  conditions  on $A$ and $B$ are then required. 

 Let $\mathcal{W}$ be the set defined by \begin{equation} \mathcal{W}= \lbrace x \in \mathcal{H}: \hspace*{0.1cm} BS(t)x=0,\hspace*{0.1cm} \forall t\geq0\rbrace. \end{equation}
 The set $\mathcal{W}$   is a closed subspace vector of $\mathcal{H}$ and  represents all unobservable states of  the pair $(A,B)$.  Let $P$ be the projection of the state space $\mathcal{H}$   on $\mathcal{W}^{\perp}$ (the orthogonal space of $\mathcal{W}$) along $\mathcal{W}$, so for all $y\in \mathcal{H}$, we can set  $y=y_{1}+y_{2}$ with $Py=y_{1}$ and $(I-P)y=y_{2}$. \\
Observing that $\mathcal{W}\subset Ker(B)$ and using the fact that $B$ is a self-adjoint operator, we deduce that  
  \begin{equation}\label{2.2}
  B\mathcal{W}\subset \mathcal{W}  \hspace*{0.5cm}\textit{and } \hspace*{0.5cm}B\mathcal{W}^{\perp} \subset \mathcal{W}^{\perp}.
  \end{equation}
 So that the projection $P$  commutes with $B$. Furthermore, we can verify that $\mathcal{W}$  is invariant under $S(t)$. Let us consider the   following assumption:\\
 \\
\textbf{$\textbf{(}\mathcal{\textbf{H}}_{\textbf{1}}\textbf{)}$}:  The set $\mathcal{W}^{\perp}$ is invariant for the semigroup $S(t)$, for all $t \geq 0$.\\
\\
Note that the  condition \textbf{$\textbf{(}\mathcal{\textbf{H}}_{\textbf{1}}\textbf{)}$} is verified if, for example,  the infinitesimal generator $A$ is self-adjoint or skew-adjoint, and  also in  the special case where  $\mathcal{W}$= $\lbrace 0 \rbrace$. Moreover, if the  hypothesis \textbf{$\textbf{(}\mathcal{\textbf{H}}_{\textbf{1}}\textbf{)}$} is  verified, then $P$ and $I-P$ commute with $A$. \\  
 According to the above discussion, we can see that  under  the assumption  \textbf{$\textbf{(}\mathcal{\textbf{H}}_{\textbf{1}}\textbf{)}$}, the initial system  (\ref{systS}) can be decomposed into the following two subsystems:
\begin{eqnarray}\label{systS1}
	\begin{cases}
		\frac{d }{d t}y_{1}(t) =A y_{1}(t)+ u(t) B y_{1}(t), \quad   \hspace*{0.1cm} t>0,\\
		y_{1}(0)=y_{01}\in  \mathcal{W}^{\perp}, \\
	\end{cases}
\end{eqnarray}
 and 
  \begin{eqnarray}\label{systS2}
	\begin{cases}
		\frac{d }{d t} y_{2}(t)=A y_{2}(t), \quad  \hspace*{0.1cm} t>0,\\
		y_{2}(0)=y_{02} \in \mathcal{W}, \\
	\end{cases}
\end{eqnarray}
where $y_{0}=y_{01}+y_{02}$, and let us  denote by $S_{1}(t)$  and $S_{2}(t)$ the restrictions of the semigroup $S(t)$ on  $\mathcal{W}^{\perp}$ and   $\mathcal{W}$, respectively.\\
  \begin{remark} 1)   By the decomposition of the state space given by \cite{Amaliki}, it is necessary to assume that $\ker(B)$ is stable by $S(t)$. On the other hand, our decomposition guarantees that the space $\mathcal{W}$ is always stable by $S(t)$. This makes our state space decomposition more suitable.\\
  	2)	It is worth noting that $S_{2}(t)$ is a $C_{0}$-semigroup on $\mathcal{W}$ whose generator  is the operator $A_{\mathcal{W}}$ induced by $A$ on $\mathcal{W}$. Moreover, under the assumption \textbf{$\textbf{(}\mathcal{\textbf{H}}_{\textbf{1}}\textbf{)}$}, $S_{1}(t)$ is a $C_{0}$-semigroup on $\mathcal{W^{\perp}}$ generated by   the operator $A_{\mathcal{W^{\perp}}}$.
  \end{remark} 
 \subsection{Well-posedness and  Stabilization results}
 Before starting our results, let us recall the definition of finite-time stabilization
of system (\ref{systS}).
\begin{definition} ([\cite{233388}, p. 112])
Let $D\subset \mathcal{H}$. The system (\ref{systS}) is finite-time  stable in $D$  if  it is Lyapunov stable and  there exists a feedback control ensuring that, for any initial state $y_{0}\in D $, a finite-time  $T=T(y_{0})$  exist   such that  $y(t)=0$ for all $t\geq T(y_{0})$. Moreover, if $D=\mathcal{H}$, then the system (\ref{systS}) is globally finite-time  stable.\\
In this case, $T^{*}(y_{0}):=\inf \hspace*{0.1cm}\{\ T(y_{0}) \geq 0: y(t)=0, \forall t \geq T(y_{0}) \} $ is called the settling-time function of the
finite-time stable system.

\end{definition} 
We show that stabilising the initial system   (\ref{systS}) in finite time, turns out to stabilizing the subsystem   (\ref{systS1}), provided that  (\ref{systS2}) satisfies some properties.\\
 Let us consider the following assumption:\\ 
 \\
\textbf{$\textbf{(}\mathcal{\textbf{H}}_{\textbf{2}}\textbf{)}$} \hspace*{0.2cm} $\langle A y, By\rangle \leq \varphi (y)\hspace*{0.1cm} \| By\|^2,  \hspace*{0.3cm} \textit{for all } y\in \mathcal{W}^{\perp}\cap D(A)$, where $ \varphi:\mathcal{W}^{\perp}\cap D(A)\rightarrow \mathbb{R}^{+}$  is such that  the mapping $\Upsilon:\mathcal{W}^{\perp}\cap D(A)\rightarrow \mathcal{H}$  defined by $\Upsilon(y)= \varphi(y)B(y)$ can be extended to $\mathcal{W}^{\perp}$  on a  globally Lipschitz function.\\ 
\\
 Furthermore,  {for such function} $\varphi$, we consider the following nonlinear feedback control:
 \begin{equation}\label{systSbnbbbbkamalfenzan}
u(t)= -  \left(\langle BPy(t), Py(t)\rangle^{-\mu}+ \varphi ( Py(t))\right) \hspace*{0.1cm}\textbf{1}_{E_{P}} (y(t)),
\end{equation}
where  $0<\mu<\frac{1}{2}$, $E_{P}=\lbrace y\in \mathcal{H},\hspace*{0.1cm}  BPy\neq 0\rbrace:={\ker(BP)}^{c}$.\\
 With the control   (\ref{systSbnbbbbkamalfenzan}), the system   (\ref{systS}) becomes
\begin{eqnarray}\label{systSresule}
	\begin{cases}
		\frac{d }{d t} y(t)=A y(t)- F(y(t)) , \quad    \hspace*{0.1cm} t>0,\\
		y(0)=y_{0},    \\
	\end{cases}
\end{eqnarray}
where  $$
F(y)= \begin{cases}  \langle BPy, Py\rangle^{-\mu}By+ \varphi(Py) By, & BPy \neq 0, \\ 0, & BPy=0.\end{cases}
$$
We begin with the following result which   introduces the existence and uniqueness of the weak  solution of the
closed-loop system  (\ref{systSresule}).
\begin{theorem}\label{theo1}  Let $A$ generate  a quasi-contractive semigroup $ S(t)$ on $\mathcal{H}$ with  a type $\omega$ (i.e. $\langle A x, x\rangle \leq \omega \| x\|^2$, $\forall x\in D(A)$) and let   $B$ be a self-adjoint positive operator  such that assumptions \textbf{$\textbf{(}\mathcal{\textbf{H}}_{\textbf{1}}\textbf{)}$} and  \textbf{$\textbf{(}\mathcal{\textbf{H}}_{\textbf{2}}\textbf{)}$} hold. Then,  for all $y_{0}\in \mathcal{H}$,  the system (\ref{systSresule}) has a unique weak solution $y\in C([0,+\infty), \mathcal{H})$.
\end{theorem}
\textbf{Proof}.
It is clear  that the unique solution of  (\ref{systS2})  is given by  
$y_{2}(t)= S_{2}(t)y_{02} $, for all  $t\geqslant 0$. Moreover, the system  (\ref{systS1}) in closed-loop   can be written as follows
 \begin{eqnarray}\label{systS12b}
	\begin{cases}
		\frac{d}{d t} y_{1}(t) +\mathcal{A}y_{1}(t)+\Upsilon(y_{1}(t))=0,\quad  \hspace*{0.1cm} t>0,\\
		y_{1}(0)=y_{01},   \\
	\end{cases}
\end{eqnarray}
where $\mathcal{A}:D(A)\cap\mathcal{W^{\perp}}\rightarrow\mathcal{W}^{\perp}$ is defined by 
\begin{equation}\label{systS12bbv}
\mathcal{A}z= -Az +\Theta(z),
\end{equation}
and $\Theta$ is defined for all $z \in  \mathcal{W^{\perp}}$ by
$$
\Theta(z)= \begin{cases}  \langle  Bz, z\rangle^{-\mu} Bz, & Bz \neq 0, \\ 0, & Bz=0.\end{cases}
$$
To establish the well-posedness of  (\ref{systS12b}), we'll show that the operator $\mathcal{A}+\omega I$  is maximal monotone.
Since $A$   generates a  quasi-contractive semigroup, the operator $-A+\omega I$ is maximal monotone.  Furthermore, from \cite{sob} the mapping $\Theta$ is maximal monotone, thus according to   Theorem \ref{theork2}, the  operator
$\mathcal{A}+\omega I$  is maximal monotone as well.\\
Using the assumption \textbf{$\textbf{(}\mathcal{\textbf{H}}_{\textbf{2}}\textbf{)}$}, we deduce from  Theorem \ref{theork1}, the operator  $\mathcal{A}+\Upsilon+ (\omega+\kappa) I$ is  maximal monotone on $D(A)\cap \mathcal{W^{\perp}}$, where $\kappa$ is the Lipschitz constant of $\Upsilon$. According to (\cite{sobb}, p. 105),  we  conclude that for all $y_{01} \in \mathcal{W^{\perp}}$,   the system  (\ref{systS12b}) has a unique weak solution $y_{1}\in C([0,+\infty), \mathcal{W^{\perp}})$.\\
Now,  let $z=y_{1}+y_{2}$, where $y_{1}$ is the solution of  (\ref{systS12b}) and $y_{2}$ the solution of  (\ref{systS2}). Through the sum of systems (\ref{systS2})    and  (\ref{systS12b}), we have 
\begin{eqnarray}\label{systSbn}
	\begin{cases}
		\frac{d}{d t} z(t) =A z(t)+ u(t)Bz(t), \quad   \hspace*{0.1cm} t>0,\\
		z(0)=  y_{0}=y_{01}+y_{02}.\\
	\end{cases}
\end{eqnarray}
 Using the uniqueness of the decomposition of z according to the sum $\mathcal{H}=\mathcal{W}\oplus\mathcal{W^{\perp}}$,  we deduce that  for all $y_{0}\in \mathcal{H}$, the solution  $z:\mathbb{R}^{+}\rightarrow \mathcal{H}$  of  (\ref{systSbn}) is the unique  solution of  the system  (\ref{systSresule}).\\
\\
We will now introduce our first main result which concerns the global finite-time  stabilisation of  (\ref{systS}). To do this,  we begin with the following proposition which provides a necessary condition for finite-time stability of (\ref{systS}).
\begin{proposition}\label{RE1n} If the system  (\ref{systS}) is globally finite-time stable, then the semigroup $S_{2}(t) $ induced by $S(t) $ on the subspace $\mathcal{W}$  is  nilpotent.
\end{proposition}
\textbf{Proof}.
Assume that (\ref{systS}) is stable  in  a finite time $T\geq 0$, that is there is a feedback low $u(t)=f(y(t))$ for which, $y(t)=0$, for all $t\geq T$.  Let $z=y(0)\in \mathcal{W}$, then regardless the stabilizing  control $u(t)$, we have 
$y(t)=S_{2}(t)z=0$ for all $t\geq T$. Hence, the semigroup $S_{2}(t) $   is  nilpotent  on $\mathcal{W}$.\\
\begin{remark}  \label{RE1} If the semigroup $S(t)$ is one to one for some $t>0$ (as   $\mathcal{H}= \mathbb{R}^{n}$,  heat and   wave equations)  and $\mathcal{W}\neq \lbrace 0 \rbrace$, then,  for any feedback control $u(t)$, the system (\ref{systS}) cannot be  finite-time  stable on $\mathcal{W}\setminus \{0\}$.
\end{remark}  
Here, we consider the following assumption: \\

\textbf{$\textbf{(}\mathcal{\textbf{H}}_{\textbf{3}}\textbf{)}$} There exists $\gamma>0$ such that  $\gamma  \langle Bx, x\rangle \leq \lVert Bx\rVert^{2}$,     for all $x\in \mathcal{W}^{\perp}$.\\
\\
In addition, based on Proposition \ref{RE1n},  we introduce the following hypothesis:\\
\\
\textbf{$\textbf{(}\mathcal{\textbf{H}}_{\textbf{4}}\textbf{)}$} $S_{2}(t) $ is a nilpotent  semigroup (i.e there exists $\delta\geq 0$ such that $S_{2}(t)= 0$, for all $t \geq \delta$).\\
\begin{remark}
 The assumption \textbf{$\textbf{(}\mathcal{\textbf{H}}_{\textbf{3}}\textbf{)}$} is verified (in particular) if the control operator $B$ is a  projection and take $\gamma=1$.\\
\end{remark}
We can state our main result on finite-time stabilization of system (\ref{systSresule}) as follows.
\begin{theorem}\label{theo2} Assume that the hypotheses of Theorem  \ref{theo1} are satisfied
and suppose  that \textbf{$\textbf{(}\mathcal{\textbf{H}}_{\textbf{3}}\textbf{)}$} and \textbf{$\textbf{(}\mathcal{\textbf{H}}_{\textbf{4}}\textbf{)}$}  hold. Then,  the system  (\ref{systSresule}) is globally finite-time stable. \\
Furthermore, the settling time  admits the following  estimate:\\
$$
T^{*}\leq max\hspace*{0.1cm}(\frac{\langle By_{01}, y_{01}\rangle^{\mu}}{2\gamma \mu},\delta).
$$
\end{theorem}
\textbf{Proof}.
According to Theorem \ref{theo1}, the system   (\ref{systS12b}) admits a unique weak solution  for all $y_{01}\in \mathcal{W^{\perp}}$, hence, there is  a sequence $y_{01}^n$ in $D(A)\cap \mathcal{W^{\perp}}$ such that $y_{01}^n \rightarrow y_{01}$, and the system  (\ref{systS12b}) with $y_{01}^n$ as initial state possesses a strong solution $y_{1}^n \in C([0,+\infty), \mathcal{W^{\perp}})$ verifying $y_{1}^n \rightarrow y_{1}$ uniformly on $[0, T]$ (for  any $T>0$ ) as $n \rightarrow \infty$.\\
 Then,  for a.e $t>0$,  we have $y_{1}^n(t) \in D(A)\cap \mathcal{W^{\perp}}$  and   
$$ 
\begin{aligned}
\frac{\mathrm{d}}{\mathrm{d} t}y_{1}^n(t)-Ay_{1}^n(t)+ (\langle By_{1}^n(t), y_{1}^n(t)\rangle^{-\mu} +  \varphi(y_{1}^n(t))By_{1}^n(t)\textbf{1}_{E} =0,
\end{aligned}
$$
where we have used that  $Py_{1}=y_{1}$ for all $y_{1}\in \mathcal{W^{\perp}}$ and $E=\lbrace t\geq0,\hspace*{0.1cm} By_{1}^n(t)\neq 0\rbrace$.\\
Now, in order to  prove the finite-time stability of the system   (\ref{systS12b}), we consider the Lyapunov function candidate  defined by $V(y_{1})=\langle By_{1}, y_{1}\rangle $, for all $y_{1}\in \mathcal{W^{\perp}}$. The time derivative of $V$ along the trajectories  of the system  (\ref{systS12b}), yields
$$
\begin{aligned}
\frac{1}{2} \frac{d}{d t}\langle By_{1}^n(t), y_{1}^n(t)\rangle&=\langle By_{1}^n(t), Ay_{1}^n(t)+u(t)B(y_{1}^n(t))\rangle. \\
\end{aligned}
$$
By assumption \textbf{$\textbf{(}\mathcal{\textbf{H}}_{\textbf{3}}\textbf{)}$}, we find
$$
\begin{aligned}
\frac{1}{2} \frac{d}{d t}V(y_{1}^n(t))\leq 
\begin{cases}
		 -\gamma V(y_{1}^n(t))^{1-\mu}, &  By_{1}^n(t) \neq 0, \\
		0, &  By_{1}^n(t)=0.
	\end{cases}
\\
\end{aligned}
$$
Hence,  according to the comparison principle (Lemma \ref{4nb}), we have
$$ 
	\begin{cases}
	\begin{aligned}
		& V(y_{1}^n(t))^{\mu}\leq V(y_{01}^n)^{\mu} -  2\gamma\mu  t, & \text{ if }   t \leq \frac{V(y_{01}^n)^{\mu}}{ 2\gamma\mu},\\
			& By_{1}^n(t)=0, & \text{ if }  t \geq  \frac{V(y_{01}^n)^{\mu}}{2\gamma \mu}.
		\end{aligned}
	\end{cases} 
$$
Then, letting  $n\rightarrow +\infty$,  we deduce that
$$ 
	\begin{cases}
	\begin{aligned}
		& V(y_{1}(t))^{\mu}\leq V(y_{01})^{\mu} -  2\gamma\mu  t,  \text{ if }   t \leq \frac{V(y_{01})^{\mu}}{ 2\gamma\mu},\\
			& By_{1}(t)=0,  \text{ if }  t \geq T_{1}= \frac{V(y_{01})^{\mu}}{2\gamma \mu}.
		\end{aligned}
	\end{cases} 
$$
Thus, the system (\ref{systS12b}) reads as follows for  $t\geq T_{1}$ 
$$
	\begin{cases}
		\frac{d }{d t} y_{1}(t)=A y_{1}(t), \quad   \hspace*{0.1cm} t\geq T_{1},\\
		y_{1}(T_{1})=  y_{1}(T_{1}^{-}), \\
	\end{cases}
$$
which implies that, $y_{1}(t)=S_{1}(t-T_{1})y_{1}(T_{1})$, for all $t\geq T_{1}$. Then,  we deduce that   $y_{1}(T_{1})\in \mathcal{W}$ and hence  $y_{1}(t)=0$,  for all $t\geq T_{1}$. This  proves  the global finite-time stability of  (\ref{systS12b}). Moreover, by the assumption \textbf{$\textbf{(}\mathcal{\textbf{H}}_{\textbf{4}}\textbf{)}$}, there exists $\delta\geq 0$ such that $y_{2}(t)= 0$, for all $t \geq \delta$. Consequently,  $y(t)=y_{1}(t)+y_{2}(t)= 0$, for all $t\geq max\hspace*{0.1cm}(T_{1},\delta)$.\\
 The operator $A$ generate  a quasi-contractive semigroup, then 
\begin{equation}\label{CON} 
\begin{aligned}
 \frac{d }{d t}\left\|y(t)\right\|^{2} &\leq  2  \hspace*{0.1cm}\omega \left\|y(t)\right\|^{2},
\end{aligned}
\end{equation}
from which, it comes
\begin{equation}\label{CONN}
\left\|y(t)\right\| \leq \left\|y_{0}\right\| \hspace*{0.1cm}\operatorname{e}^{\frac{\omega \varrho}{2}}, \forall t \in [0, \varrho),
\end{equation}
with $\varrho = max\hspace*{0.1cm}(T_{1},\delta)$. Hence (\ref{systSresule}) is Lyaponov stable. Consequently,  the system  (\ref{systSresule}) is globally finite-time stable at a settling time $T^{*}\leq \varrho = max\hspace*{0.1cm}(T_{1},\delta$).\\
\\
From  the previous results, we can deduce the following finite-time stability result:
\begin{corollary}\label{KoaaazzbVV}
Let the bounded operator  $B$ be  self-adjoint  and positive, and let assumptions \textbf{$\textbf{(}\mathcal{\textbf{H}}_{\textbf{1}}\textbf{)}$}, \textbf{$\textbf{(}\mathcal{\textbf{H}}_{\textbf{2}}\textbf{)}$} and \textbf{$\textbf{(}\mathcal{\textbf{H}}_{\textbf{3}}\textbf{)}$} hold. Then, the system  (\ref{systSresule}) with initial state  $y_{0} \in  \mathcal{W^{\perp}}$ vanishes at a settling time $T^{*}\leq\frac{\langle By_{0}, y_{0}\rangle^{\mu}}{2\gamma\mu}$.
\end{corollary}
\begin{remark}\label{Remark}
 1) In the case  $\mathcal{W}$= $\lbrace 0 \rbrace$, the output finite-time stability of some classes of abstract bilinear systems has been studied in \cite{sob}.\\
2) If the system  (\ref{systS}) is stable in finite-time, then under the assumption \textbf{$\textbf{(}\mathcal{\textbf{H}}_{\textbf{1}}\textbf{)}$} both the systems  (\ref{systS1}) and  (\ref{systS2}) are finite-time stable.\\  
 3) If the function $\varphi$ is constant, then the condition \textbf{$\textbf{(}\mathcal{\textbf{H}}_{\textbf{2}}\textbf{)}$} is obviously verified. Moreover, if it is zero, then  the system  (\ref{systS}) is  finite-time stable under the following feedback control:
\begin{equation}\label{systSbnbbbbkamalfenza}
u(t)= -    \langle BPy(t), Py(t)\rangle^{-\mu} \hspace*{0.1cm}\textbf{1}_{E_{P}} (y(t)), 
\end{equation}
where $E_{P}=\ker(BP)^{c}$ and $0<\mu<\frac{1}{2}$.\\
\end{remark}
The following proposition discusses the global  finite-time stability  of any (eventually multiple)  weak solution of system  (\ref{systS})    under the following feedback control: 
\begin{equation}\label{systSbnbbbbkamal}
 u(t)= -  (\langle BPy(t), Py(t)\rangle^{-\mu}+ \frac{\langle A Py(t), BPy(t)\rangle}{\lVert BPy(t)\rVert^{2}}) \hspace*{0.1cm}\textbf{1}_{E_{P}}(y(t)),
\end{equation}
where $E_{P}=\ker(BP)^{c}$ and $0<\mu<1$.
\begin{proposition}\label{Koaaazzbbb}
Let assumptions \textbf{$\textbf{(}\mathcal{\textbf{H}}_{\textbf{3}}\textbf{)}$} and \textbf{$\textbf{(}\mathcal{\textbf{H}}_{\textbf{4}}\textbf{)}$} hold  and assume that $B$ is  a self adjoint  positive operator. Then,  every weak solution of the system in closed loop   (\ref{systS})-(\ref{systSbnbbbbkamal}), vanishes  at a finite time.
In addition, the settling time  admits the following  estimate:\\
$$
T^{*}\leq \frac{\langle By_{01}, y_{01}\rangle^{\mu}}{2\gamma\mu}+\delta.
$$
\end{proposition}
\textbf{Proof}.
For all $y_{0}\in \mathcal{H}$, we write  $y_{0}=y_{01}+y_{02}$ with $y_{01}\in \mathcal{W^{\perp}}$  and $y_{02}\in \mathcal{W}$. Using  (\ref{2.2}),  the  weak solution $y(t)$ of  (\ref{systS}) under the state feedback  (\ref{systSbnbbbbkamal}) can be uniquely written as $y(t)=y_{1}(t)+y_{2}(t)$, where  $y_{1}(t)$ is the weak  solution of the following system:
\begin{eqnarray}\label{systS11}
	\begin{cases}
		\frac{d }{d t}y_{1}(t) =PA y_{1}(t)+ u(t) B y_{1}(t), \quad    \hspace*{0.1cm} t>0,\\
		y_{1}(0)=y_{01}\in  \mathcal{W}^{\perp}, \\
	\end{cases}
\end{eqnarray}
and $y_{2}(t)$ is the weak  solution of  the  system
\begin{eqnarray}\label{systS22}
	\begin{cases}
		\frac{d }{d t} y_{2}(t)=A y_{2}(t)+(I-P)A y_{1}(t), \quad   \hspace*{0.1cm} t>0,\\
		y_{2}(0)=y_{02} \in \mathcal{W}. \\
	\end{cases}
\end{eqnarray}
 According to  (\ref{2.2}), we get that  $\langle Ay_{1}, By_{1}\rangle=\langle PAy_{1}, By_{1}\rangle$, for all $y_{1}\in \mathcal{W}^{\perp}\cap D(A)$. Hence,  we deduce from the proof of Theorem \ref{theo2} that system  (\ref{systS11}) is globally finite-time stable and the settling time can be estimated as  $T_{1}\leq\frac{\langle By_{01}, y_{01}\rangle^{\mu}}{2\gamma\mu}$.\\
We have the following variation of constants formula, for all $t\geq s \geq 0$:
$$
y_{2}(t)= S_{2}(t-s)y_{2}(s)+ \int_{s}^{t}   S_{2}(t-\tau) (I-P)Ay_{1}(\tau)d\tau.
$$
It follows from the condition \textbf{$\textbf{(}\mathcal{\textbf{H}}_{\textbf{4}}\textbf{)}$} that there exists $\delta\geq 0$ such that $S_{2}(t)= 0$, for all $t \geq \delta$.\\
Consequently, for $s=\frac{\langle By_{01}, y_{01}\rangle^{\mu}}{2\gamma\mu}$, we have 
$y_{2}(t)=0$, for all $t\geq \frac{\langle By_{01}, y_{01}\rangle^{\mu}}{2\gamma\mu}+\delta$. We conclude that the system  (\ref{systS}) is globally finite-time stable at a settling
time $T^{*} \leq\frac{\langle By_{01}, y_{01}\rangle^{\mu}}{2\gamma\mu}+\delta$.

\section{Finite-time stability of linear systems}
In  this section,  we will use the results obtained for bilinear systems to study the problem of  finite-time  stability  of the following linear system:
 \begin{eqnarray}\label{systSb1B}
	\begin{cases}
		\frac{d }{d t} y(t)=A y(t)+Lv(t), \quad    \hspace*{0.1cm}  t>0,\\
		y(0)=y_{0},   \\
	\end{cases}
\end{eqnarray}
where $A:D(A)\subset \mathcal{H}\rightarrow  \mathcal{H}$ is the infinitesimal generator of a $C_{0}$-semigroup $S(t)$ on $ \mathcal{H}$ and  the control input  $v$ is a  $\mathcal{U}$-valued function, with $\mathcal{U}$ being a real Hilbert space. Moreover, $L:\mathcal{U}\rightarrow \mathcal{H}$ is a bounded linear operator.\\
 Then, the space $\mathcal{W}$ is  given by 
$$ \mathcal{W}= \lbrace x \in \mathcal{H}, \hspace*{0.1cm} L^{*}S(t)x=0,\hspace*{0.1cm} \forall t\geq0\rbrace,$$
where $L^{*}$  is the adjoint operator of $L$.\\
The following result provides the  finite-time stability of  (\ref{systSb1B}):
\begin{theorem}\label{KoaaazzbVVbbbbbbbbbbvvb}
 Let  $S(t)$ be a quasi-contractive semigroup on $\mathcal{H}$ that satisfies  the assumptions \textbf{$\textbf{(}\mathcal{\textbf{H}}_{\textbf{1}}\textbf{)}$} and \textbf{$\textbf{(}\mathcal{\textbf{H}}_{\textbf{4}}\textbf{)}$} and suppose   that hypotheses \textbf{$\textbf{(}\mathcal{\textbf{H}}_{\textbf{2}}\textbf{)}$} and \textbf{$\textbf{(}\mathcal{\textbf{H}}_{\textbf{3}}\textbf{)}$} hold for  $B=LL^{*}$. Then, the system  (\ref{systSb1B}) is  globally finite-time stable by the following state feedback law:
\begin{equation}\label{systSbnbbbbbb}
v(t)= -  \left(\frac{L^{*}Py(t)}{\| L^{*}Py(t)\|^{2\mu}} + \varphi (Py(t))L^{*}Py(t) \right) \hspace*{0.1cm}\textbf{1}_{E_{P}}  (y(t)),
\end{equation} 
where $E_{P}=\lbrace y\in \mathcal{H},\hspace*{0.1cm}  L^{*}Py\neq 0\rbrace$ and $0<\mu<\frac{1}{2}$. Furthermore, the settling time  admits the following  estimate:\\
$$
T^{*}\leq max ( \delta, \frac{\| L^{*}y_{01}\|^{2\mu}}{2\gamma \mu}).
$$
\end{theorem}
\textbf{Proof}. 
Under the control  (\ref{systSbnbbbbbb}), the system  (\ref{systSb1B}) has the form of  (\ref{systS}) if we take $B=LL^{*}$ and $v(t)=u(t) L^{*}y(t)$.\\
It is clear that the operator $B$ is positive and self-adjoint from $\mathcal{H}$ to $\mathcal{H}$. \\Applying the result of Theorem \ref{theo1} and  Theorem \ref{theo2} to system  (\ref{systSb1B}), we obtain that there exists a   time  $ T=max ( \delta, \frac{\| L^{*}y_{01}\|^{2\mu}}{2\gamma \mu}) $, for which we have   $y(t)=0$, for all $t\geq T$, where $\delta\geq 0$ is such that   $S_{2}(t)=0$, for all $t\geq \delta$.\\

In the following, we provide another control to stabilise the system   (\ref{systSb1B}), for which the assumptions  on the operator $B$  are systematically satisfied, and  the system   (\ref{systSb1B}) in   closed-loop admits a unique weak solution without assuming  that \textbf{$\textbf{(}\mathcal{\textbf{H}}_{\textbf{2}}\textbf{)}$} holds.\\
Let us  introduce the control
\begin{equation}\label{inssafk}
v(t)=-\langle Py(t), \zeta\rangle|\langle Py(t), \zeta\rangle|^{-2 \mu} \varpi-\frac{\left\langle Py(t), A^* \zeta\right\rangle}{\|\zeta\|^2} \varpi,
\end{equation}
where $\zeta \in D\left(A^*\right)$ and $\varpi \in U$ to be determined, and  $P$  is the projection on  $\mathcal{W}^{\perp}=P\mathcal{H}$ along $\mathcal{W}=(I-P)\mathcal{H}$. \\
Let us consider the operator $A_\zeta=A-\frac{\left\langle . \hspace*{0.1cm},\hspace*{0.1cm}P A^* \zeta\right\rangle}{\|\zeta\|^2} \zeta$ with domain $D\left(A_\zeta\right)=D(A)$. Observing that the linear operator $N y=\frac{\left\langle y, PA^* \zeta\right\rangle}{\|\zeta\|^2} \zeta$ is bounded, we deduce that $A_\zeta$ generates a quasi-contractive
semigroup  $T(t)$   which is given  by: 
$$
\begin{aligned}
T(t)y&= S(t)y-\frac{1}{\|\zeta\|^2}\int_{0}^{t}S(t-\tau) \left\langle T(\tau)y , PA^* \zeta\right\rangle \zeta d\tau,  \\
\end{aligned}
$$  
for all $y\in \mathcal{H}$ and all $t\geq 0$. \\
Here,  let $\mathcal{W}$ be the set defined by $\mathcal{W}=\left\{y \in \mathcal{H} ;\langle T(t) y, \zeta\rangle=0, \forall t \geq\right.$ $0\}$. Then, we have the following   finite-time stability result:
\begin{proposition}\label{pro2}
Suppose that\\
1) The set $\mathcal{W}$ is invariant for the semigroup $T^*(t)$, for all $t \geq 0$, and   the assumption \textbf{$\textbf{(}\mathcal{\textbf{H}}_{\textbf{4}}\textbf{)}$}  holds,\\
2) there exist a non-zero $\zeta \in D\left(A^*\right)$ and $\varpi \in U$ such that $L \varpi=\zeta$. \\
Then, for all $y_{0}\in  \mathcal{H}$,  the system  (\ref{systSb1B}) with the state feedback  (\ref{inssafk}) is  finite-time stable and the settling time admits the estimate $T^{*}\leq max\lbrace \frac{\left|\left\langle y_{01}, \zeta\right\rangle\right|^{ \mu}}{\mu\|\zeta\|^2},\delta\rbrace $.
\end{proposition}
\textbf{Proof}.
\textit{\underline{  Step 1 :}} Well-posedness\\
  By decomposing the system  (\ref{systSb1B}) with the state feedback  (\ref{inssafk}) and according to assumption \textbf{$\textbf{(}\mathcal{\textbf{H}}_{\textbf{1}}\textbf{)}$} and  the fact that $\zeta\in \mathcal{W^{\perp}}$, we get
 \begin{equation}\label{systS2b1}
\begin{aligned}
\begin{cases}
		\frac{d }{d t}y_{1}(t)=A y_{1}(t)-\langle y_{1}(t), \zeta\rangle|\langle y_{1}(t), \zeta\rangle|^{-2 \mu} \zeta-\frac{\left\langle y_{1}(t), A^* \zeta\right\rangle}{\|\zeta\|^2} \zeta, \\
y_{1}(0)=y_{01} \in \mathcal{W}^{\perp},
	\end{cases}
\end{aligned}
\end{equation}
and 
\begin{eqnarray}\label{systS2b}
	\begin{cases}
		\frac{d }{d t}y_{2}(t)=A y_{2}(t), \\
y_{2}(0)=y_{02} \in \mathcal{W}.
	\end{cases}
\end{eqnarray}
Then, for all $y_{02} \in \mathcal{W}$,  the unique weak solution of  (\ref{systS2b}) is given by $y_{2}(t)=S(t)y_{02}$, for all $t\geq 0$.\\
For the well-posedness of system  (\ref{systS2b1}),  we'll use the fact that the  mapping  $\mathcal{B} y_{1}=\langle y_{1}, \zeta\rangle|\langle y_{1}, \zeta\rangle|^{-2 \mu} \zeta$  is maximal monotone on $\mathcal{H}$ (see \cite{sob}).\\
 We have  that $A_{\zeta}$   generates a   quasi-contractive
semigroup $T(t)$ with  type $\omega+ \|N\|$,  ($\omega$ is the type of $S(t)$), then  from Theorem \ref{theork2}  it follows that $\mathcal{A}+(\omega+ \|N\|)I$, with $\mathcal{A}= \mathcal{B}- A_\zeta$,  is a maximal monotone operator .\\
 Then, the system  (\ref{systS2b1})  has a unique weak solution  $y_{1}\in C([0,+\infty), \mathcal{W^{\perp}})$, for all $y_{01}\in\mathcal{W^{\perp}}$, (see \cite{sobb}, p. 105). \\
 We conclude that $y(t)=y_{1}(t)+y_{2}(t)$ is the unique weak solution of  (\ref{systSb1B}).\\
\textit{\underline{  Step 2 :}}   finite-time stability\\
To study the finite-time stability of  (\ref{systS2b1}), let's consider the following  Lyaponuv function: $V(y_{1})= \langle y_{1}, \zeta\rangle ^{2}$, for all $ y_{1}\in \mathcal{W^{\perp}}$. \\Similarly to the proof of Theorem \ref{theo2} we can show that
$$
	\begin{cases}
		\left|\left\langle  y_{1}(t), \zeta\right\rangle\right|^{2 \mu} \leq\left|\left\langle y_{01}, \zeta\right\rangle\right|^{2 \mu}-2\mu\|\zeta\|^2 t,0 \leq t \leq \frac{\left|\left\langle y_{01}, \zeta\right\rangle\right|^{2 \mu}}{2\mu\|\zeta\|^2},\\
			\left|\left\langle y_{1}(t), \zeta\right\rangle\right|=0, \hspace*{0.1cm} t \geq T_{1}=\frac{\left|\left\langle y_{01}, \zeta\right\rangle\right|^{2 \mu}}{2\mu\|\zeta\|^2}.\\
	\end{cases}
$$
Consequently, the system  (\ref{systS2b1}) can be written as 
\begin{eqnarray}\label{systS2bKAMAL}
	\begin{cases}
		\frac{d }{d t}y_{1}(t)= A_\zeta y_{1}(t), \quad   \hspace*{0.1cm} t\geq T_{1},\\
		y_{1}(T_{1})=  y^{*}. \\
	\end{cases}
\end{eqnarray}
Thus the solution of (\ref{systS2b1})  is given by  $y_{1}(t)= T(t-T_{1})y^{*}$, for all $t\geq T_{1}$. Hence, $y^{*} \in \mathcal{W} \cap \mathcal{W^{\perp}}$. Therefore, we have  $y_{1}(t)= 0$, for all $t\geq T_{1}$. Hence,  with  assumption \textbf{$\textbf{(}\mathcal{\textbf{H}}_{\textbf{4}}\textbf{)}$}, the  unique solution of  (\ref{systSb1B})  with the control (\ref{inssafk}) is  given by  
$y(t)= 0$, for all  $t\geq max\lbrace T_{1},\delta\rbrace$. This completes the proof of the proposition.
\begin{remark}\label{hg}
If $y_{0}\in \mathcal{W}^{\perp}$, we can establish the finite-time stability of the system   (\ref{systSb1B}) with the  feedback control  (\ref{inssafk}) on $\mathcal{W}^{\perp}$ by applying Proposition  \ref{pro2}, even without the assumption  \textbf{$\textbf{(}\mathcal{\textbf{H}}_{\textbf{4}}\textbf{)}$}.
\end{remark}

\section{Applications}

In this part, we illustrate our results by several examples of parabolic and hyperbolic described by bilinear and linear control systems. 
\subsection{Heat equation}
 Let us consider the one-dimensional heat equation
\begin{equation}\label{systSbbbbbbbbbbv}
\begin{cases}\begin{aligned}
		& \frac{\partial y(t)}{\partial t} =\Delta  y(t)+ u(t)By(t),   &t\in\mathbb{R}^{+}\\
		& y(t,0)=y(t,1)=0,  & t\in \mathbb{R}^{+}\\
		& y(0, x)=y_{0}, & x \in (0,1),
		\end{aligned}
	\end{cases} 
 \end{equation}
 where $y (t)=y (t,\cdot) $ is the state. The system  (\ref{systSbbbbbbbbbbv}) has been studied for a positive definite control operator  by Polyakov et al (see \cite{17}) and  Ouzahra (see \cite{2333}). Our aim here is to show that when the control operator  $ B$ is positive not necessary  definite and not coercive, this system can be stabilized in finite time  on  a subspace of the state space $\mathcal{H}=L^{2}(0,1)$.\\
We can write  (\ref{systSbbbbbbbbbbv}) in the form of  (\ref{systS}) on $L^{2}(0,1)$  if we set $Ay=\Delta y$ for $y\in D(A)=H_{0}^{1}(0,1) \cap H^{2}(0,1)$,  the spectrum of $A$ is given by the
simple eigenvalues $\lambda_{j}= -(j\pi)^{2}$ and eigenfunctions $\varphi_{j}=\sqrt 2 sin(j\pi x)$, $j\geq 1$. Furthermore $A$ generates a contraction semigroup given by $S(t) y=\sum\limits_{j=1}^{\infty} \exp \left(\lambda_{j} t\right)\langle y, \phi_j\rangle \phi_j$, for all $y\in L^{2}(0,1)$ and all $t\geq 0$.\\ Consider the bounded operator $B :L^{2}(0,1) \rightarrow L^{2}(0,1)$ defined by $By= y- \langle  y, \varphi_{1}\rangle \varphi_{1}$. Then it is clear that $B^{2}=B =B^{*}\geq 0$. Moreover, for all $y\in L^{2}(0,1)$, we have  
$$
\begin{aligned}
	&\hspace*{0.5cm}BS(t)y= 	0, \hspace*{0.3cm}\forall t \geq 0 \\
	&\Longleftrightarrow S(t)y -\langle S(t)y,\varphi_{1}\rangle \varphi_{1} =0,  \hspace*{0.3cm} t\geq 0\\
	&\Longleftrightarrow  y=\langle y,\varphi_{1}\rangle \varphi_{1},\\
	&\Longleftrightarrow  y\in span(\varphi_{1}).\\
\end{aligned}
$$
 Then, we deduce that $\mathcal{W}=span(\varphi_{1})$ and $\mathcal{W}^{\perp}=span(\varphi_{j}, j\geq 2)$ . Furthermore, since the operator $A$ is self-adjoint,  then the assumption \textbf{$\textbf{(}\mathcal{\textbf{H}}_{\textbf{1}}\textbf{)}$}  holds.  \\For all $y\in D(A)\cap \mathcal{W}^{\perp}$, we have $ \langle Ay, By\rangle = \langle Ay, y\rangle \leq 0.$  Consequently,  applying Corollary \ref{KoaaazzbVV} and point (3) of the Remark \ref{Remark}, we deduce that the system  (\ref{systSbbbbbbbbbbv}) is finite-time stable for all $y_{0} \in \mathcal{W}^{\perp}$ and  $\mu\in (0,\frac{1}{2})$, by the following feedback control:
\begin{equation}\label{567}
u(t)= \begin{cases} -\| y(t)\|^{-2\mu}, &  y(t) \neq 0 \\ 0, & y(t)=0,\end{cases}
\end{equation} 
Moreover, the semigroup $S(t)$ is one to one. As highlighted in Remark  \ref{RE1}, the system (\ref{systSbbbbbbbbbbv}) cannot be stabilized in finite time within   $\mathcal{W}=span(\varphi_{1})\setminus \{0\}$ for any alternative choice of the feedback control $u(t)$.
\subsection{Coupled transport-heat equations} 
 In this example, we  study  the global finite-time stability of the following  coupled heat-transport equations: 
\begin{equation}\label{Coupled}
\begin{cases}\frac{\partial \Phi}{\partial t}(x, y, t)=\frac{\partial^2 \Phi}{\partial x^2}(x, y, t)+\frac{\partial^2 \Phi}{\partial y^2}(x, y, t), & (x, y) \in \Omega, \\ \frac{\partial \Psi}{\partial t}(x, y, t)=-\frac{\partial \Psi}{\partial x}(x, y, t)-\frac{\partial \Psi}{\partial y}(x, y, t), & (x, y) \in \Omega, \\ \frac{\partial}{\partial \eta} \Phi(\xi, t)=0, &  \xi \in \partial \Omega  \\ \Psi(\xi, t)=0, & \xi\in \Gamma_{\vartheta},\end{cases}
\end{equation}
where  $\Omega=(0,1)^2$,  $\frac{\partial}{\partial \eta}$ is the  normal derivative on $\partial \Omega$ and $\Gamma_{\vartheta}=\lbrace x\in \partial \Omega: \vartheta\cdot \gamma (x)<0\rbrace$, where $\gamma$ is the outward normal vector to $\partial \Omega$ and $\vartheta=(1,1)$. \\
Let's consider the state space $\mathcal{H}=L^2(\Omega) \times L^2(\Omega)$. The spectrum of the Laplacian operator in $\Omega$ with Neumann boundary conditions are given by the simple eigenvalues\\ $\lambda_{j, k}= -\left(j^2+k^2\right) \pi^2$, $(j, k)\in \mathbb{N}^{2} $ corresponding to eigenfunctions $\varphi_{j, k}(x, y)=2 \cos (j \pi x)\times  \cos (k \pi x)$ for all $(j, k)\in \mathbb{N}^{2}  {\backslash\lbrace{(0,0)}\rbrace}$ and $\varphi_{0, 0}(x, y)=1$.\\
 Moreover, the Laplacian operator $\Delta$ with domain $\mathcal{D}_1=\lbrace\Phi \in H^2(\Omega): \frac{\partial }{\partial \eta}\Phi=0$ on $\partial \Omega\rbrace$ generates on $L^2(\Omega)$ the following contraction semigroup:  $
\mathcal{T}(t)\Phi=\sum\limits_{j,k=0}^{+\infty} e^{\lambda_{j, k} t}\left\langle\Phi, \varphi_{j, k}\right\rangle \varphi_{j, k} 
$, for all $\Phi\in L^2(\Omega)$. Let $\nabla$  design the gradient operator, then the operator $\vartheta \cdot \nabla$ with domain $\mathcal{D}_2=\left\{\psi \in W^{1,2}(\Omega)\right.: \left.\left.\psi(\xi)\right|_{\Gamma_{\vartheta}}=0\right\}$ generates in $L^2(0,1)$ the left translation semigroup defined by
$$\mathcal{S}(t) \psi(\xi)= \begin{cases}\psi(\xi+v t), & \text { if } \xi+v t \in \Omega \\ 0, & \text { else }\end{cases}$$
 which is a contraction semigroup.\\
Then, the operator  $A=\left(\begin{array}{cc}\Delta & 0 \\ 0 & v \cdot \nabla\end{array}\right)$ with domain $D(A)=\mathcal{D}_1 \times \mathcal{D}_2$, generates on $\mathcal{H}$ the following contraction semigroup $\mathcal{H}(t)=\left(\begin{array}{cc}\mathcal{T}(t) & 0 \\ 0 & \mathcal{S}(t)\end{array}\right)$.\\
Let $h >0$  and $\omega=(0,h)^2 $.  Clearly,  the system (\ref{Coupled}) is unstable.  Here, we will prove the finite-time stabilization of   the following bilinear control coupled system:\\
 \begin{equation}\label{Coupled1}
\begin{cases}\frac{\partial \Phi}{\partial t}(x, y, t)=\frac{\partial^2 \Phi}{\partial x^2}(x, y, t)+\frac{\partial^2 \Phi}{\partial y^2}(x, y, t)+ u(t) \Phi(x, y, t), & (x, y) \in \Omega,\\ \frac{\partial \Psi}{\partial t}(x, y, t)=-\frac{\partial \Psi}{\partial x}(x, y, t)-\frac{\partial \Psi}{\partial y}(x, y, t)+ u(t) \chi_{\omega}\Psi(x, y, t), & (x, y) \in \Omega, \\ \frac{\partial}{\partial \eta} \Phi(\xi, t)=0, &  \xi \in \partial \Omega,  \\ \Psi(\xi, t)=0, & \xi\in \Gamma_{\vartheta},\end{cases}
\end{equation}
where  $\chi_{\omega}$ indicates the characteristic function of $\omega$. This system takes the form (\ref{systS}) if we  consider $B=\left(\begin{array}{cc}I  & 0 \\0 & \chi_{\omega}\end{array}\right)$, which is 
  a self-adjoint and  positive operator in $\mathcal{H} $  and satisfies  the assumption \textbf{$\textbf{(}\mathcal{\textbf{H}}_{\textbf{3}}\textbf{)}$}.   Then,   the space $\mathcal{W}$ is given by  
$$
\mathcal{W}= \left\lbrace (0,\Psi) \in \mathcal{H}: \Psi|_{\omega}=0\right\rbrace,
$$
and so,
 $$
\mathcal{W^{\perp}}= \left\lbrace (\Phi,\Psi) \in \mathcal{H}: \Psi|_{\omega^{c}}=0  \right\rbrace.
$$
 For all $z= (\Phi,\Psi) \in\mathcal{W^{\perp}}$, we have  that $\mathcal{H}(t)z\in \mathcal{W^{\perp}}$, thus, the assumption \textbf{$\textbf{(}\mathcal{\textbf{H}}_{\textbf{1}}\textbf{)}$} is true. In addition,  the assumption \textbf{$\textbf{(}\mathcal{\textbf{H}}_{\textbf{2}}\textbf{)}$} is satisfied with $\varphi=0$. \\
Let $Y=(0, \phi)\in \mathcal{W}$, then $\mathcal{H}(t)Y= (0, \mathcal{S}(t)\phi)= (0,0)$, for all $t\geq 1$. Hence, the assumption \textbf{$\textbf{(}\mathcal{\textbf{H}}_{\textbf{4}}\textbf{)}$} is fulfilled.\\
Consequently,  applying the results of Theorem \ref{theo2},   the feedback control
$$
u(t)= - \left| \int_{\Omega} \Phi^2(x, y) dx dy+\int_{\omega} \Psi^2(x, y) dx dy\right|^{-\mu} \hspace*{0.1cm}\textbf{1}_{E}(\Phi,\chi_{\omega}\Psi),
$$
with  $\mu \in (0,\frac{1}{2})$ and $E=\lbrace (\Phi,\Psi)\in \mathcal{H}: (\Phi,\Psi)\neq (0,0)\rbrace$, guarantees the global finite-time stability of the system (\ref{Coupled1}).

\subsection{ Wave equation}
 Let $\Omega=(0,1)$, $Q=\Omega \times( 0,+\infty)$ and $q\geq 1$. Let us consider the one-dimensional wave equation:
\begin{eqnarray}\label{systSbbb}
	\begin{cases}
		\frac{\partial^{2} }{\partial t^{2}}y(x,t) =\Delta  y(x,t)+ u(t) \sum\limits_{i=1}^{q} \langle \varphi_{i}(.),\frac{\partial }{\partial t}y(.,t) \rangle \varphi_{i}(x), ~~ &   (x,t)\in Q,\\
		y(0,t)=y(1,t)=0, ~~ &  t\in \mathbb{R}^{+}, \\
		y(x, t)=y_{0}, ~~ & x \in \Omega.
	\end{cases}
\end{eqnarray}  
 Here where $y (t)=y (t,\cdot) $ is the state and denote by $\mathcal{H}$ the product Hilbert space $H_0^1(\Omega) \times L^2(\Omega)$ with the inner product defined by:
$$
\langle\left(y_1, y_2\right),\left(z_1, z_2\right)\rangle= \langle y_1, z_1\rangle_{H_0^1(\Omega)}+\langle y_2, z_2\rangle_{L^2(\Omega)}, 
$$
for all $(y_{i}, z_{i})_{i=1,2} \in \mathcal{H}$.\\
 We can write  (\ref{systSbbb}) in the first-order differential version  (\ref{systS}) on $\mathcal{H}$ as 
\begin{eqnarray}\label{systSZES}
	\begin{cases}
		\frac{\partial}{\partial t} \Phi(t) =A \Phi(t)+ Lv(t) , \quad    t>0,\\
		\Phi(0)=\Phi_{0},   \\
	\end{cases}
	\end{eqnarray}
where $\begin{aligned}
\Phi(t) & =\left(\begin{array}{c}
y(t) \\
\frac{\partial y(t)}{\partial t}
\end{array}\right) 
\end{aligned}$ and $
\mathcal{A} =\begin{pmatrix}  0&I \\ A&0 \end{pmatrix}$,  where $A=\Delta  =\partial^2 / \partial x^2$, with $D(A)=H_0^1(0,1) \cap H^2(0,1)$.  
The control operator $B=\left(\begin{array}{ll} 0 & 0 \\ 0 & G(.) \end{array}\right)$, where  $G: L^2(\Omega) \rightarrow L^2(\Omega)$ is defined by $G (y)= \sum\limits_{i=1}^{q}\langle   \varphi_{i},y \rangle \varphi_{i}$.  Then, the control operator $B$ is self-adjoint, positive in $\mathcal{H}$ and the assumption  \textbf{$\textbf{(}\mathcal{\textbf{H}}_{\textbf{3}}\textbf{)}$} is verified.\\
Moreover, the spectrum of the operator $-A$ with Dirichlet conditions is given by the simple eigenvalues $\lambda_j=(j \pi)^2$  corresponding to the eigenfunctions  $\varphi_j(x)=\sqrt{2} \sin (j \pi x), \forall j \in \mathbb{N}^*$. \\
Let $y=\left(y_1, y_2\right) \in \mathcal{H}$, with $y_1=\sum\limits_{j=1}^{\infty} \alpha_j \varphi_j$ and $y_2=\sum\limits_{j=1}^{\infty} \lambda_j^{\frac{1}{2}} \beta_j \varphi_j$, where $\left(\alpha_j, \beta_j\right) \in \mathbb{R}^2$, $ \forall j \geq 1$. Separation of variables yields
$$
S(s) y=\sum_{j=1}^{\infty}\left(\begin{array}{c}
\alpha_j \cos \left(\lambda_j^{\frac{1}{2}} s\right)+\beta_j \sin \left(\lambda_j^{\frac{1}{2}} s\right) \\
-\alpha_j \lambda_j^{\frac{1}{2}} \sin \left(\lambda_j^{\frac{1}{2}} s\right)+\beta_j \lambda_j^{\frac{1}{2}} \cos \left(\lambda_j^{\frac{1}{2}} s\right)
\end{array}\right) \varphi_j. 
$$
 Moreover, for every  $y\in \mathcal{W}$, we have   
$$
\begin{aligned}
&\hspace*{1cm} BS(t)y= 	0,\hspace*{0.4cm}  \forall t\geq 0, \\
&\Rightarrow  -\alpha_i \lambda_i^{\frac{1}{2}} \sin (\lambda_i^{\frac{1}{2}} t)+\beta_i \lambda_i^{\frac{1}{2}} \cos (\lambda_i^{\frac{1}{2}} t)  =0,\; \; \forall t\geq 0,\;\; \forall  i \in \lbrace 1,..., q \rbrace,\\
&\Rightarrow \alpha_i= \beta_i=0, \hspace*{1cm}   \forall  i \in \lbrace 1,..., q \rbrace. \\
\end{aligned}
$$
It follows that 
$$
 \mathcal{W}=  \lbrace y= (y_1, y_2)\in \mathcal{H},  y_1= \sum_{j\geq q+1} \alpha_j \varphi_j \hspace*{0,3cm} \textit{and} \hspace*{0,3cm} y_2= \sum_{j\geq q+1} \beta_j \varphi_j, \hspace*{0,3cm} \left(\alpha_j, \beta_j\right) \in \mathbb{R}^2 \rbrace.
$$
From which, we can  easily deduce that
$$
 \mathcal{W^{\perp}}=  \lbrace y= (y_1, y_2)\in \mathcal{H},  y_1= \sum_{j=1}^{q} \alpha_j \varphi_j \hspace*{0,3cm} \textit{and}\hspace*{0.4cm} y_2=  \sum_{j=1}^{q} \beta_j \varphi_j, \hspace*{0,3cm} \left(\alpha_j, \beta_j\right) \in \mathbb{R}^2 \rbrace.
$$ 
Since the operator $\mathcal{A}$   is skew-adjoint,  the assumption \textbf{$\textbf{(}\mathcal{\textbf{H}}_{\textbf{1}}\textbf{)}$} is fulfilled.\\
Moreover, for any $y=(y_1, y_2)\in D(A) \cap \mathcal{W^{\perp}}$ we have 
$$
\begin{aligned}
\left \langle \mathcal{A} y, By \right \rangle &= \left\langle (y_2, Ay_1), (0, \sum_{i=1}^{q} \langle   \varphi_{i},y_2 \rangle \varphi_{i})\right\rangle\\	
&= \sum_{i=1}^{q} \langle Ay_{1},  \varphi_{i}\rangle \langle   \varphi_{i},y_2 \rangle\\
&=  \sum_{i=1}^{q} -\lambda_i\alpha_i \beta_i. \\
\end{aligned}
$$
Therefore, the hypothesis \textbf{$\textbf{(}\mathcal{\textbf{H}}_{\textbf{2}}\textbf{)}$}  is satisfied with $\varphi(y)=K=\max\left\lbrace  |\frac{ \alpha_i}{\beta_i}|, \beta_i \neq 0,  i \in \lbrace 1,..., q \rbrace\right\rbrace$, for all $y\in D(A) \cap \mathcal{W^{\perp}}$.
Then, using Corollary \ref{KoaaazzbVV}, the feedback control
$$
u(t)=- ((\sum\limits_{i=1}^{q}\langle   \varphi_{i},\frac{\partial}{\partial t}y(.,t) \rangle ^2)^{- \mu}+ K)  \hspace*{0.1cm}\textbf{1}_{\ker(B)^{c}} (y(.,t)),
$$
where $0<\mu<\frac{1}{2}$,  guarantees the  finite-time stability of the system (\ref{systSbbb})  for all $y_{0}\in \mathcal{W^{\perp}}$.\\
Observing that  the semigroup $S(t)$ is one to one. Then, applying   Remark \ref{RE1},   we conclude  that  (\ref{systSbbb}) is cannot be  stable in finite-time  for all $y_{0}\in \mathcal{W}\setminus \{0\}$ for any
 choice of the feedback control $u(t)$.
\subsection{ Beam equation}

 The beam equation is a second-order hyperbolic equation that generally describe a set of physical problems such as oscillating systems like wave equations and vibrating beams (see \cite{29}).  Strong and uniform  stabilization of beam equation  was investigated in \cite{L} via nonlinear feedback, and  the exponential stabilisation result  has been obtained  in \cite{Ouzahra2010}. Here, we will use  the Proposition \ref{pro2} to prove the finite-time stability of a linear  beam equation. \\
  Let $\Omega \subset \mathbb{R}^n$, $n\geq 1$, denote a bounded open domain with sufficiently smooth boundary $\Gamma$ and consider the  following system:
\begin{eqnarray}\label{systSZESvvvvvbbbbbb}
\begin{cases}
\frac{\partial^2 y(x, t)}{\partial t^2}=-\frac{\partial^4 y(x, t)}{\partial x^4}+v (t)h, & (x, t) \in \Omega \times(0,+\infty), \\ y(x, t)=  \frac{\partial^2 y(x, t)}{\partial x^2}=0, & (x, t) \in \Gamma \times(0,+\infty), \\ y(x, 0)=y_0, & x \in \Omega.\end{cases}
\end{eqnarray}
Let 
$$
F=\frac{\partial^4}{\partial x^4},   \text { with } D(F)=\left\{ y \in L^2(\Omega) \,\middle|\, \frac{\partial^4 y}{\partial x^4} \in L^2(\Omega),\; y(x) = \frac{\partial^2 y(x) }{\partial x^2}= 0,\; x\in \Gamma   \right\} \text { and } h\in D(F)\setminus \{0\}.
$$
 The spectrum of the operator $F$ is given by the simple eigenvalues $\lambda_j=(j \pi)^4$,  corresponding to the eigenfunctions  $\varphi_j(x)=\sqrt{2} \sin (j \pi x), \forall j \in \mathbb{N}^*$.\\
This system can be written   in the form   (\ref{systS}) if we set\\
 $\bullet$ $\mathcal{H} =\left(H^2(\Omega) \cap H_0^1(\Omega)\right) \times L^2(\Omega)$,\\
  $\bullet$ $
A=\left(\begin{array}{cc}
0 & I \\
-F& 0
\end{array}\right) 
$ with $D(A)= D(F) \times D(F^{\frac{1}{2}})$.\\
 $\bullet$ $L: \mathbb{R} \rightarrow \mathcal{H}$ is defined by $L \lambda=\lambda\left(\begin{array}{l}0 \\ h\end{array}\right)$.\\
 $\bullet$  $h\in \ker(F)$.\\
 Then,  the operator $A$ is skew-adjoint. Then the assumption \textbf{$\textbf{(}\mathcal{\textbf{H}}_{\textbf{1}}\textbf{)}$}  holds. Moreover,  the adjoint operator $L^*: \mathcal{H} \rightarrow \mathbb{R}$ of $L$ is given by $L^*\left(\begin{array}{l}y_1 \\ y_2\end{array}\right)=\left\langle y_2, h\right\rangle_{L^2(\Omega)}$. Then,  we have  $L L^*=\left(\begin{array}{ll}0 & 0 \\ 0 & \Upsilon(.)\end{array}\right)$, where  $\Upsilon: L^2(\Omega) \rightarrow L^2(\Omega)$ is defined by $\Upsilon (y)=$ $\langle y, h\rangle_{L^2(\Omega)} h$.\\
Here, the state space  is equipped  with the following inner product:
 $$
\begin{aligned} 
 \left\langle\left(\begin{array}{l}y_1 \\ y_2\end{array}\right),\left(\begin{array}{l}z_1 \\ z_2\end{array}\right)\right\rangle=\left\langle F^{1 / 2}y_1, F^{1 / 2}z_1\right\rangle_{L^2(\Omega)}+
 \left\langle y_2, z_2\right\rangle_{L^2(\Omega)}, \hspace*{0.4cm} \forall (y_{i}, z_{i})_{i=1,2} \in \mathcal{H}.\\
 \end{aligned}$$
The  output finite-time stability of the system  (\ref{systSZESvvvvvbbbbbb}) with respect to the output operator  $C: \mathcal{H}\rightarrow  \mathbb{R}^{2}$ defined by $C\left(\begin{array}{l}y_1 \\ y_2\end{array}\right)= \left(\begin{array}{l}\langle y_{1}, h\rangle \\ \langle y_{2}, h\rangle\end{array}\right)$ has been studied in \cite{sob}. Here, we will established that the trajectory of (\ref{systSZESvvvvvbbbbbb}) vanishes at a finite time.\\ 
 We have
$$
\begin{aligned}
\mathcal{W}= \lbrace y=(y_{1},y_{2})\in \left(H^2(\Omega) \cap H_0^1(\Omega)\right) \times L^2(\Omega) : 
\langle y_{1}, h\rangle  =\langle y_{2}, h\rangle=0 \rbrace.\\
\end{aligned}
$$
Consequently, applying the result of Proposition \ref{pro2} and Remark \ref{hg}, we deduce that the following feedback control law:
\begin{equation}\label{567}
v(t)= \begin{cases} -\left\langle \frac{\partial}{\partial t}y(.,t),  h\right\rangle_{L^2(\Omega)}\left|\left\langle \frac{\partial}{\partial t}y(.,t), h\right\rangle_{L^2(\Omega)}\right|^{-2 \mu}, &  \langle \frac{\partial}{\partial t}y(.,t),  h\rangle_{L^2(\Omega)} \neq 0 \\ 0, &  \langle \frac{\partial}{\partial t}y(.,t),  h\rangle_{L^2(\Omega)}=0,\end{cases}
\end{equation} 
where $0<\mu<\frac{1}{2}$, ensures  the finite-time stability of  (\ref{systSZESvvvvvbbbbbb}) on $\mathcal{W}^{\perp}$.\\
As noted  in Remark  \ref{RE1}, the system (\ref{systSbbbbbbbbbbv}) cannot be stabilized in finite time within $\mathcal{W}\setminus \{0\}$ for any  choice of the feedback control $v(t)$.

\section{Conclusion}
% In this paper, we have tackled the problem of global finite-time stability for abstract bilinear and linear systems. While existing research has primarily focused on systems that are weakly observable, our study broadens this scope by addressing bilinear systems that may not be observable and have control operators that are positive but not necessarily definite. We explicitly derive the stabilizing feedback control for such systems, leveraging a novel criterion based on state-space decomposition. The effectiveness of this approach is further demonstrated by applying the results to both parabolic and hyperbolic equations, highlighting the method’s broad applicability and potential impact. This study  concerns   the case of bounded systems. However, the modelling   may give raise to unbounded control operators. This is the case when the action is exercised on a part of the boundary or concentrated in a point of the geometrical domain of the system at hand. Thus,   the question of finite-time stability of such systems is of great interest.

In this paper, we addressed the problem of global finite-time stability for abstract bilinear and linear systems. While much of the existing literature focuses on systems that are weakly observable, our study broadens this perspective by considering bilinear systems that may not be observable and whose control operators are positive but not necessarily definite and not necessarily corercive. We provided an explicit construction of stabilizing feedback controls based on a novel state-space decomposition criterion. The effectiveness and versatility of this method were demonstrated through its application to both parabolic and hyperbolic partial differential equations, underscoring its wide-ranging potential. Our current analysis has been limited to the case of systems with bounded control operators. However, many practical hyperbolic systems involve unbounded control actions, such as those exerted at the boundary or localized at a point within the spatial domain. This introduces new mathematical challenges, particularly in ensuring finite-time stability under such control configurations. Investigating these scenarios remains an important direction for future work. Moreover, another promising line of research involves extending the current framework to time-delay systems, which frequently arise in networked control and distributed parameter systems. Developing robust finite-time stabilizing controllers that can handle input and state delays while maintaining the guarantees provided by our approach would be of significant practical and theoretical interest
\bibliographystyle{plain}
\bibliography{Bib}
\end{document}